\documentclass [11pt]{article}
\usepackage{mathtools,amsfonts,amssymb,amsthm, bm, nccmath}
\usepackage{epsfig}
\usepackage{graphicx}
\allowdisplaybreaks
\usepackage{enumerate}
\usepackage{cite}
\begin{document}
	\thispagestyle{empty}
	\null\vspace{-1cm}
	\medskip
	\vspace{1.75cm}
	\centerline{\Large\textbf{{Generalization of Cauchy type theorems for matrix Polynomials}}}
	~~~~~~~~~~~~~~~~~~~~~~~~~~~~~~~~~~~~~~~~~~~~~~~~~~~~~~~~~~~~~~~~~~~~~~~~~~~~~~~~~~~~~~~~~~~~~~~~~~~~~~~~~~~~~~~~~~~~~~~~~~~~~~~~~~~~~~~~~~~~~~~~~~~~~~~~~~~

	\centerline{\bf {Idrees Qasim}}
\centerline {Department of Mathematics, National Institute of Technology, Srinagar, India-190006}
\centerline {idreesf3@nitsri.ac.in}

	~~~~~~~~~~~~~~~~~~~~~~~~~~~~~~~~~~~~~~~~~~~~~~~~~~~~~~~~~~~~~~~~~~~~~~~~~~~~~~~~~~~~~~~~~~~~~~~~~~~~~~~~~~~~~~~~~~~~~~~~~~~~~~~~~~~~~~~~~~~~~~~~~~~~~~~~~
	\vskip0.1in
	\noindent \textbf{Abstract}: In this paper, we find bounds for the eigenvalues of matrix polynomials. In particular, we find generalizations of Cauchy's classical Theorem for distribution of eigenvalues of matrix polynomial.

	~~~~~~~~~~~~~~~~~~~~~~~~~~~~~~~~~~~~~~~~~~~~~~~~~~~~~~~~~~~~~~~~~~~~~~~~~~~~~~~~~~~~~~~~~~~~~~~~~~~~~~~~~~~~~~~~~~~~~~~~~~~~~~~~~~~~~~~~~~~~~~~~~~~~~~~~~~
	
	\noindent {{\bf Keywords:} Complex Polynomials, Matrix polynomial, Zeros, Eigenvalues, Cauchy Theorem.
	\vspace{0.15in}
	~~~~~~~~~~~~~~~~~~~~~~~~~~~~~~~~~~~~~~~~~~~~~~~~~~~~~~~~~~~~~~~~~~~~~~~~~~~~~~~~~~~~~~~~~~~~~~~~~~~~~~~~~~~~~~~~~~~~~~~~~~~~~~~~~~~~~~~~~~~~~~~~~~~~~~~~~
	
	\noindent {{\bf Mathematics Subject Classification (2020):} } 12D10, 15A18, 30C15\\
	\vspace{0.05in}
	~~~~~~~~~~~~~~~~~~~~~~~~~~~~~~~~~~~~~~~~~~~~~~~~~~~~~~~~~~~~~~~~~~~~~~~~~~~~~~~~~~~~~~~~~~~~~~~~~~~~~~~~~~~~~~~~~~~~~~~~~~~~~~~~~~~~~~~~~~~~~~~~~~~~~~~~~
	\section{Introduction} 
	Let $p(z):=\sum_{j=0}^{m}a_jz^j$ be a polynomial of degree $m$ with complex coefficients. One of the basic theorems of mathematics is the Fundamental Theorem of Algebra,
	according to which, ``\textit{every polynomial of degree $n\ge 1$ has exactly $n$ zeros in the complex
	plane}". This theorem does not however say anything regarding the location of zeros of
	a polynomial. Finding zeros of a polynomial with complex coefficients presents a significant challenge, not confined to mathematics alone but extending across various scientific, engineering, and technological domains. However, according to the Abel-Ruffini Theorem (for further details, refer to \cite{AR}), a universally applicable method for finding the roots of a polynomial of degree 5 or higher solely through arithmetic operations and radicals does not exist. Consequently, it becomes imperative to ascertain the domains where these roots are situated. In this pursuit, Cauchy's classical theorem, first articulated in 1829 (referenced as \cite{CAU}), emerges as a valuable instrument for analyzing the distribution of polynomial roots. The theorem states:\\
	\textbf{Theorem A.} All the zeros of the polynomial $p(z)=a_nz^n+a_{n-1}z^{n-1}+\cdots+a_1z+a_0$, $a_n\neq 0$, with complex coefficients lie in
	$$\{z;|z|<s\}\subset \{z;|z|<1+M\},$$
	where $M=\max\limits_{1\le j\le n-1}\left|\dfrac{a_j}{a_n}\right|$ and $s$ is the unique positive solution of the equation
	$$|a_n|z^n-|a_{n-1}|z^{n-1}-\cdots-|a_1|z-|a_0|=0.$$
	
	The problem of locating some or all the zeros of a given polynomial as
	a function of its coefficients is of long standing interest in mathematics. This fact can
	be deduced by glancing at the references in the comprehensive books of Marden \cite{LM}
	and Milovanovi´c, Mitrinovi´c and Rassias \cite{MMR},  Rahman and Schmeisser \cite{RS} and by noting the abundance of recent	publications on the subject.
\\
	
	\noindent Consider $\mathbb{M}_{n\times n}$ as the set comprising all $n\times n$ matrices with entries from $\mathbb{C}$. We define a matrix polynomial, denoted as $P:\mathbb{C}\rightarrow \mathbb{M}_{n\times n}$, as a function expressed by
	$$P(z):=\sum_{j=0}^{m}A_jz^j,~~~A_j\in  \mathbb{M}_{n\times n}.$$
	If $A_m\neq 0$, then $P(z)$ is termed a matrix polynomial of degree $m$. An eigenvalue of $P(z)$, denoted by $\lambda$, is a value for which there exists a non-zero vector $u\in\mathbb{C}^n$ such that $P(\lambda)u=0$. In such a case, $u$ is called an eigenvector of $P(z)$. The problem of finding a number $\lambda\in \mathbb{C}$ and a non-zero vector $u\in \mathbb{C}^n$ satisfying $P(\lambda)u=0$ is referred to as a Polynomial Eigenvalue Problem ($PEP$). When $m=1$, this reduces to a Generalized Eigenvalue Problem ($GEP$) represented as $Au=\lambda Bu$. Furthermore, if $B$ equals the identity matrix $I$, it becomes the standard eigenvalue problem $Au=\lambda u$. To see more about matrix polynomials see \cite{GLR}.\\
	
	The computation of eigenvalues of a matrix polynomial poses a significant challenge, often leading to approximate solutions due to numerical errors. Iterative methods are commonly employed to tackle this challenge (for references, see \cite{11}). When computing pseudospectra of matrix polynomials, it becomes crucial to delineate a region within the complex plane encompassing the eigenvalues of interest. Bounds can aid in delineating such a region (for references, see \cite{7}). Moreover, bounding the eigenvalues can expedite numerical computations by aiding in selecting appropriate numerical methods and reducing the number of iterations required in various numerical algorithms.
	
	Note that, if $A_0$ is singular then 0 is an eigenvalue of $P(z)$, and if $A_m$ is singular then 0 is an eigenvalue of the matrix polynomial $z^mP(1/z)$. Therefore, to locate the eigenvalues of these matrix polynomials, we always assume	that $A_0$ and $A_m$ are non-singular.
	
	\noindent \textbf{Notation}. Throughout this paper, $||\cdot||$ denotes a subordinate matrix norm.\\

	The following extensions of Theorem A to matrix polynomial were obtained in \cite{CAM,7, 2}.
	
	\noindent \textbf{Theorem B.} Let $P(z):=\sum_{j=0}^{m}A_jz^j$, $\det(A_m)\neq 0$ be a matrix polynomial. Then the eigenvalues of $P(z)$ lie in $|z|\le \rho$, where $\rho$ is a unique positive root of the equation
	
	$$||A_m^{-1}||^{-1}z^m-||A_{m-1}||z^{m-1}-\cdots-||A_1||-||A_0||=0.$$
	
	\noindent \textbf{Theorem C.} Let $P(z):=\sum_{j=0}^{m}A_jz^j$, $\det(A_m)\neq 0$ be a matrix polynomial. Then the eigenvalues of $P(z)$ lie in $|z|<1+M$, where $M=||A_m^{-1}||\max\limits_{0\le j\le m-1}||A_j||.$\\
	
	In this paper, we find some new bounds for the eigenvalues of matrix polynomials.

	\section{Main Results}
	We first prove the following result:\\
\textbf{Theorem 1.} All the eigenvalues of the matrix polynomial $P(z):=\sum\limits_{j=0}^{m}A_jz^j$ lie in
\begin{equation}
	|z|<\left[\frac{1}{2}\left\{1+(1+4\alpha^q_p)^{1/2}\right\}\right]^{1/q},\label{eq:MN1}
\end{equation}
where
$$\alpha_p:=\left\{\sum_{r=0}^{m}\left(\frac{||A_{m-1}A_{m-r}-A_mA_{m-r-1}||}{||A^{-2}_m||^{-1}}\right)^p\right\}^{1/p},$$
$p>1$, $\dfrac{1}{p}+\dfrac{1}{q}=1$, $A^-2=(A^2)^{-1}$ and $A_{-1}=0$.
\begin{proof}
	We consider
\begin{align*}
	P^*(z)&=(A_{m-1}-A_mz)P(z)\\
	&=\left(A_{m-1}-A_mz\right)\left(A_mz^m+A_{m-1}z^{m-1}+\cdots+A_1z+A_0\right)\\
	&=-A^2_mz^{m+1}+\left(A^2_{m-1}-A_mA_{m-2}\right)z^{m-1}+\cdots+A_{m-1}A_0,
\end{align*}
therefore, for a unit eigenvector $u$, we have
\begin{align*}
	||P^*(z)u||&=\left|\left|\left\{-A^2_mz^{m+1}+\left(A^2_{m-1}-A_mA_{m-2}\right)z^{m-1}+\cdots+A_{m-1}A_0\right\}u\right|\right|\\
	&\ge \left|\left|(A^2_mz^{m+1})u\right|\right|-\left|\left|\sum_{r=1}^{m}(A_{m-1}A_{m-r}-A_mA_{m-r-1})z^{m-r}\right|\right|\\
	&\ge ||A^{-2}_m||^{-1}|z|^{m+1}-\sum_{r=1}^{m}||A_{m-1}A_{m-r}-A_mA_{m-r-1}|||z|^{m-r}\\
	&=||A_m^{-2}||^{-1}|z|^{m+1}\left[1-\sum_{r=1}^{m}\frac{||A_{m-1}A_{m-r}-A_mA_{m-r-1}||}{||A^{-2}_m||^{-1}}.\frac{1}{|z|^{r+1}}\right].
\end{align*}	
By Holder's inequality,
\begin{align*}
	\sum_{r=1}^{m}\frac{||A_{m-1}A_{m-r}-A_mA_{m-r-1}||}{||A^{-2}_m||^{-1}}.\frac{1}{|z|^{r+1}}&\le \left\{\sum\limits_{r=1}^{m}\left(\frac{||A_{m-1}A_{m-r}-A_mA_{m-r-1}||}{||A^{-2}_m||^{-1}}\right)^p\right\}^{1/p}\\
	&~~~~~~~~~~~~~~~~~~~~~~~~~~~~~\times\left\{\sum_{r=1}^{m}\frac{1}{|z|^{q(r+1)}}\right\}^{1/q}\\
	&=\alpha_p\left\{\sum_{r=1}^{m}\frac{1}{|z|^{q(r+1)}}\right\}^{1/q}.
\end{align*}
Hence for $|z|>1$,
\begin{align*}
	||P^*(z)u||&\ge ||A^{-2}_m||^{-1}|z|^{m+1}\left\{1-\frac{\alpha_p}{\left\{|z|^q(|z|^q-1)\right\}^{1/q}}\right\}\\
	&>0,
\end{align*}
if 
$$|z|^{2q}-|z|^q-\alpha^q_p\ge 0.$$
That is if,
$$|z|\ge \left[\frac{1}{2}\left\{1+(1+4\alpha^q_p)^{1/2}\right\}\right]^{1/q}.$$
Therefore, those eigenvalues of $P^*(z)$ and hence of $P(z)$ lie in
$$|z|<\left[\frac{1}{2}\left\{1+(1+4\alpha^q_p)^{1/2}\right\}\right]^{1/q}$$
whose modulus is greater than 1. But, those eigenvalues of $P(z)$ whose modulus is less than or equal to 1 already lie in the region \eqref{eq:MN1}. This proves the theorem.
\end{proof}
\noindent \textbf{Remark 1}. A result of Mohammad \cite{Mohammad} for the distribution of zeros of complex polynomial is a special case of Theorem 1 when we take $n=1$.\\

Next, we prove the following result:\\
\noindent \textbf{Theorem 2.} If $P(z):=\sum\limits_{j=0}^{m}A_jz^j$ is a matrix polynomial of degree $m$, then for any $p$ and $q$ such that $p>1$, $q>1$ with $\dfrac{1}{p}+\dfrac{1}{q}=1$, all the eigenvalues of $P(z)$ lie in
\begin{equation}
	|z|<\left[1+\left(A_p\right)^q\right]^{1/q},\label{eq:N01}
\end{equation}                          
where $A_p=\left[\sum\limits_{j=0}^{m-1}\left(\dfrac{||A_j||}{||A_m^{-1}||^{-1}}\right)^p\right]^{1/p}.$
\begin{proof}
Let $u$ be a unit vector, then we have
\begin{align*}
	||P(z)u||&\ge ||(A_mz^m)u||-\left|\left|\sum_{j=0}^{m-1}(A_{j}z^j)u\right|\right|\\
	&\ge||A_m^{-1}||^{-1}|z|^m-\sum_{j=0}^{m-1}||A_j|||z|^j\\
	&=||A_m^{-1}||^{-1}|z|^m\left[1-\sum_{j=0}^{m-1}\frac{||A_j||}{||A_m^{-1}||^{-1}}.\frac{1}{|z|^{m-j}}\right].\\	
\end{align*}
Using Holder's inequality,
\begin{equation}
	||P(z)u||\ge ||A_m^{-1}||^{-1}|z|^m\left[1-\left\{\sum_{j=0}^{m-1}\left(\frac{||A_j||}{||A_m^{-1}||^{-1}}\right)^p\right\}^{1/p}.\left\{\sum_{j=0}^{m-1}\frac{1}{|z|^{(m-j)q}}\right\}^{1/q}\right].\label{eq:0002}	
\end{equation}
If $|z|>1$, then
\begin{align*}
	\sum_{j=0}^{m-1}\frac{1}{|z|^{(m-j)q}}&=\sum_{j=1}^{m}\frac{1}{|z|^{jq}}\\
	&<\sum_{j=1}^{\infty}\frac{1}{|z|^{jq}}\\
	&=\frac{1}{|z|^q-1}.
\end{align*}
Therefore, from \eqref{eq:0002}, we have for $|z|>1$
$$||P(z)u||	>||A_m^{-1}||^{-1}|z|^m\left[1-\frac{A_p}{(|z|^q-1)^{1/q}}\right].$$
Thus for $|z|>1,$
$$||P(z)u||>0~~~\text{if}~~~|z|^q-1\ge (A_p)^q.$$
That is, if
$$|z|\ge \left[1+(A_p)^q\right]^{1/q}.$$
This shows that those eigenvalues of $P(z)$ whose modulus is greater than 1 lie in $|z|<\left[1+(A_p)^q\right]^{1/q}$. But those eigenvalues of $P(z)$ whose modulus is less or equal to one already satisfy this inequality. This proves the theorem.
\end{proof}

\noindent \textbf{Remark 2}. An easy computation gives $$\left[\left\{1+\left(\sum_{j=0}^{m-1}\left(\dfrac{||A_j||}{||A_m^{-1}||^{-1}}\right)^p\right)^{q/p}\right\}\right]^{1/q}\le (1+m^{q/p}M^q)^{1/q},$$
where $M=||A_m^{-1}||\max\limits_{0\le j\le n-1}||A_j||$.\\ 
In view of above inequality, if we let $p\rightarrow \infty$ in \eqref{eq:N01}, so that $q\rightarrow 1$ and $(1+n^{q/p}M^q)^{1/q}\rightarrow 1+M$, we get Theorem C.\\

\noindent \textbf{Remark 3}. A result of Montel \cite{Montel} for the bounds for the zeros of complex polynomials follows from Theorem 2 when we take $n=1$.\\

Next, we prove:\\
\noindent \textbf{Theorem 3.} If $P(z)=A_mz^m+A_pz^p+\cdots+A_1z+A_0$, $0\le p\le m-1$ is a matrix polynomial of degree $m$ and 
$$M=\max_{0\le j\le p}\frac{||A_j||}{||A_m^{-1}||^{-1}},$$
then all the eigenvalues of $P(z)$ lie in $|z|<k$, where $k>1$ is the unique positive root of the trinomial equation
\begin{equation}
	x^{m-p}-x^{m-p-1}-M=0.\label{eq:01}
\end{equation}
\begin{proof}
	For a unit vector $u$, we have
	$$||P(z)u||\ge ||A_m^{-1}||^{-1}|z|^m\left[1-\left(\frac{||A_p||}{||A_m^{-1}||^{-1}}\frac{1}{|z|^{m-p}}+\cdots+\frac{||A_1||}{||A_m^{-1}||^{-1}}\frac{1}{|z|^{m-1}}+\frac{||A_0||}{||A_m^{-1}||^{-1}}\frac{1}{|z|^m}\right)\right]$$
Now,
$$\frac{||A_j||}{||A_m^{-1}||^{-1}}\le M,~~~~j=0,1,2,\dots,p,$$
it follows that for $|z|>1$,
\begin{align*}
	||P(z)u||&\ge ||A_m^{-1}||^{-1}|z|^m\left[1-\frac{M}{|z|^{m-p}}\left(\sum_{j=0}^{p}\frac{1}{|z|^j}\right)\right]\\
	&>||A_m^{-1}||^{-1}|z|^m\left[1-\frac{M}{|z|^{m-p-1}}\frac{1}{|z|-1}\right]\\
	&>0,
\end{align*}
if
$$|z|^{m-p}-|z|^{m-p-1}-M\ge 0.$$
This implies $||P(z)u||>0$ if $|z|\ge k$, where $k$ is the unique positive root of the trinomial equation defined by \eqref{eq:01}. Hence all the eigenvalues of $P(z)$ whose modulus is greater than 1 lie in $|z|<k.$ Since all those eigenvalues whose modulus is less or equal to 1 already lie in $|z|<k$, we conclude that all the eigenvalues of $P(z)$ lie in $|z|<k$ and this completes the proof of the theorem.
\end{proof}
\noindent \textbf{Remark 4} For $p=m-1$, Theorem 3 reduces to Theorem C.\\

Finally, we give the following improvement of Theorem C.\\
\noindent \textbf{Theorem 4.} All the eigenvalues of the matrix polynomial $P(z):=\sum\limits_{j=0}^{m}A_jz^j$ of degree $m$ lie in the circle
\begin{equation}
	|z|<\dfrac{1}{2}\left[1+(1+4M)^{1/2}\right],\label{eq:N02}
\end{equation}where 
$$M:=\max\limits_{1\le j\le m}\dfrac{||A_{m-1}A_{m-j}-A_mA_{m-j-1}||}{||A_n^{-2}||^{-1}},~~~~A_{-1}=0,~~A^{-2}=(A^2)^{-1}.$$
\begin{proof}
	Consider
	\begin{align*}
		P^*(z)&=(A_{m-1}-A_{m}z)P(z)\\
		&=-A^2_mz^{m+1}+\sum_{j=1}^{m}(A_{m-1}A_{m-j}-A_mA_{m-j-1})z^{m-j}.
	\end{align*}
If $|z|>1$, then for a unit vector $u$, we have
\begin{align*}
	||P(z)u||&\ge ||A^2_mz^{m+1}u||-\left|\left|\sum_{j=1}^{m}(A_{m-1}A_{m-j}-A_mA_{m-j-1})z^{m-j}\right|\right|\\
	&\ge ||A_m^{-2}||^{-1}|z|^{m+1}\left[1-\sum_{j=1}^{m}\frac{||A_{m-1}A_{m-j}-A_mA_{m-j-1}||}{||A_m^{-2}||^{-1}}\frac{1}{|z|^{1+j}}\right]\\
		&>||A_m^{-2}||^{-1}|z|^{m+1}\left[1-M\sum_{j=1}^{\infty}\frac{1}{|z|^{j+1}}\right]\\
		&=||A_m^{-2}||^{-1}|z|^{m+1}\left[1-\frac{M}{|z|(|z|-1)}\right]\\
		&>0,
	\end{align*}
if 
$$|z|^2-|z|-M\ge 0.$$
That is, if
$$|z|>\dfrac{1}{2}\left[1+(1+4M)^{1/2}\right].$$
Thus the only eigenvalues of $P^*(z)$ for $|z|>1$ are those satisfying \eqref{eq:N02}. But as all the eigenvalues of $P^*(z)$ for $|z|\le 1$ already satisfy \eqref{eq:N02}, the result follows for $P^*(z)$. Since the eigenvalues of $P(z)$ are also eigenvalues of $P^*(z)$, the result follows.
\end{proof} 

\section{Conclusion}
We have generalized Cauchy type bounds for the zeros of a complex polynomial to eigenvalues of a matrix polynomials for subordinate matrix norm. We explored some new bounds for eigenvalues of matrix polynomial and extended some results of complex polynomials to matrix polynomials. The subject worth exploring involves extending results from scalar polynomials to matrix polynomials, beyond those discussed in this paper. While it's evident that not all scalar results can be extended, there's no denying that some can be extended. A valuable reference for scalar results, including some that may not be widely recognized, is \cite{RS}.
\section{Declaration}
\begin{itemize}
	\item Availability of data and materials: Not applicable
	\item Competing interests: The author declare that we have no competing interest.
	\item Funding: The author does not receive any funding for this research
	\item Authors' contributions - provide individual author contribution: Not applicable
	\item Authors' information: Department of Mathematics, National Institute of Technology, Srinagar, India, 190006. Email: idreesf3@nitsri.ac.in
\end{itemize}
	
\end{document}